\documentclass[11pt,leqno]{article}
\usepackage{latexsym}
\usepackage{amsmath,amssymb}
\usepackage{amsthm}
\usepackage{enumerate}
\title{Products of finite connected subgroups}
\author{M.~P. G\'{A}LLEGO\thanks{M. Pilar G\'{a}llego passed away on the 22nd of May, 2019. We had the privilege to work with her and to
experience her insight and generosity to share her ideas. We miss her as a collaborator and
friend.}, P. HAUCK, L. S. KAZARIN,\\
A. MART\'{I}NEZ-PASTOR
 and M.~D. P\'{E}REZ-RAMOS}
\date{}

\newcommand{\So}{\mathcal{S}}
\newtheorem{teor}{Theorem}[section]
\newtheorem{lem}[teor]{Lemma}
\newtheorem{pro}[teor]{Proposition}
\newtheorem{cor}[teor]{Corollary}
\theoremstyle{definition}
\newtheorem{rem}[teor]{Remark}
\newtheorem{de}[teor]{Definition}
\newtheorem{notation}[teor]{Notation}

\DeclareMathOperator{\Syl}{Syl}

\begin{document}
\maketitle

\begin{abstract}
For a non-empty class of groups $\cal L$, a finite group $G = AB$ is said to
be an $\cal L$-connected product of the subgroups $A$ and $B$ if $\langle a,  b\rangle \in \cal L$ for all
$a \in A$ and $b \in B$. In a previous paper, we prove that for such a product, when $\cal L = \cal S$ is the class of finite soluble groups, then
$[A,B]$ is soluble. This generalizes the theorem of Thompson which
states the solubility of finite groups whose two-generated subgroups
are soluble. In the present paper our result is applied to extend to
finite groups previous research in the soluble universe. In particular,
we characterize connected products for relevant classes of groups;
among others the class of metanilpotent groups and the
class of groups with nilpotent derived subgroup. Also we give local descriptions of relevant
subgroups of finite groups.
\medskip

\noindent
{\bf 2010 Mathematics Subject Classification.} 20D10, 20D40, 20D25. \medskip

\noindent
{\bf Keywords.} Finite groups, Products of subgroups, Two-generated subgroups, $\cal L$-connection, Fitting classes, Fitting series,  Formations.

\end{abstract}

\section{Introduction and main results}
All groups considered in this paper are assumed to be finite.
We take further previous research on the influence of two-generated subgroups on the structure of groups, in interaction with the study of products of subgroups. In \cite{HKMP} the following result is proven:

\begin{teor}\label{MT} Let the finite group $G = AB$ be the product of subgroups $A$ and $B$. Then the following statements
are equivalent:
\begin{itemize} \item[(1)] $\langle a,b \rangle$ is soluble for all $a\in A$ and $b\in B$, i.e. $A, B$ are $\So$-connected for the class $\So$ of all finite soluble groups (cf. Definition~\ref{defcon} below).
\item[(2)] For all primes $p \ne q$, all $p$-elements $a \in A$ and all $q$-elements $b \in B$, $\langle a,b \rangle$
is soluble.
\item[(3)] $[A,B] \le G_{\So}$, where $G_{\So}$ denotes the soluble radical of $G$ (i.e. the largest soluble normal subgroup of $G$).
\end{itemize}
\end{teor}

Obviously, for the special case $A=B=G$, the following well-known result of J.~Thompson is derived:

\begin{teor}\textup{(Thompson,  \cite{TH, FL95})} A finite group $G$ is soluble if and only if every two-generated subgroup of $G$ is soluble.
\end{teor}

Thompson's theorem has been generalized and sharpened in various
ways. We mention in particular the extension of R.~Guralnick, K.~Kunyavski\u{\i},
E.~Plotkin and A.~Shalev, which  describes the elements in the soluble radical $G_{\cal S}$ of a finite group $G$.

\begin{teor}\label{G_S}\textup{(Guralnick, Kunyavski\u{\i},
Plotkin, Shalev,  \cite{GKPS})} Let $G$ be a finite group and let $x\in G$. Then $x\in G_{\cal S}$ if and only if the subgroup $\langle x,y\rangle$ is soluble for all $y\in G$.
\end{teor}

Again the application of Theorem~\ref{MT}, with $A = G$, and
$B = \langle x\rangle \le G$, assures that $\langle x\rangle G_{\cal S}$ is a normal (soluble) subgroup of G under the hypothesis in statement (1), and implies Theorem~\ref{G_S}.
\smallskip

This shows how an approach involving factorized groups provides a more general setting for local-global questions related to two-generated subgroups.
A first extension of Thompson's theorem for products of groups was obtained by A.~Carocca~\cite{C}, who proved the solubility of $\cal S$-connected products of soluble subgroups. This way the following general connection property turns out to be useful:

\begin{de}\label{defcon}\textup{(Carocca, \cite{C96})} Let $\cal L$ be a non-empty class of groups. Subgroups $A$ and $B$ of a group $G$ are $\cal L$-connected if $\langle a,b\rangle\in \cal L$ for all $a\in A$ and $b\in B$. If $G=AB$ we say that $G$ is the $\cal L$-connected product of the subgroups $A$ and $B$.
\end{de}

Structure and properties of $\cal N$-connected products, for the class $\cal N$ of finite nilpotent groups, are well known (cf. \cite{Ncon1,Ncon2,Ncon3}); for instance, $G=AB$ is an $\cal N$-connected product of $A$ and $B$ if and only if $G$ modulo its hypercenter is a direct product of the images of $A$ and $B$. Apart from the above-mentioned results regarding $\cal S$-connection, corresponding studies for the classes $\mathcal N^2$ and $\cal NA$ of metanilpotent groups, and groups with nilpotent derived subgroup, respectively, have been carried out in \cite{GHP08,GHP08b}; in \cite{GHP11} connected products for the class $\mathcal S_\pi \mathcal S_\rho$ of finite soluble groups that are extensions of a normal $\pi$-subgroup by a $\rho$-subgroup, for arbitrary sets of primes $\pi$ and $\rho$, are studied. The class $\mathcal S_\pi \mathcal S_\rho$ appears in that reference as the relevant case of a large family of formations, named nilpotent-like Fitting formations, which comprise a variety of classes of groups, such as the class of $\pi$-closed soluble groups, or groups with Sylow towers with respect to total orderings of the primes. A study in \cite{GHP16} of connected subgroups, for the class of finite nilpotent groups of class at most $2$, contributes generalizations of classical results on $2$-Engel groups.

In the present paper,  as an application of Theorem~\ref{MT}, we show that main results in \cite{GHP08,GHP08b,GHP11}, proved for soluble groups, remain valid for arbitrary finite groups. In particular, we characterize connected products for some relevant classes of groups (see Theorem~\ref{Theorem}). For instance, we prove that for a finite group $G=AB$,
the subgroups $A$ and $B$ are $\mathcal N^2$-connected if and only if $A/F(G)$ and $B/F(G)$ are $\cal N$-connected, which means that for all $a\in A$ and $b\in B$, $\langle a,b\rangle ^{\cal N}\le F(\langle a,b\rangle)$ if and only if for all $a\in A$ and $b\in B$, $\langle a,b\rangle ^{\cal N}\le F(G)$, where for any group $X$, $F(X)$ denotes the Fitting subgroup of $X$, and $X^{\cal N}$ denotes the nilpotent residual of $X$, i.e. the smallest normal subgroup of $X$ with nilpotent quotient group. When we specialize our results to suitable factorizations, as mentioned above, we derive descriptions of the elements in $F_{k}(G)$, the radical of a group $G$ for the class $\mathcal N^k$ of soluble groups with nilpotent length at most $k\ge 1$, as well as the elements in the hypercenter of $G$ modulo $F_{k-1}(G)$, in the spirit of the characterization of the soluble radical in Theorem~\ref{G_S} (see Corollaries~\ref{Cor4-08},~\ref{Cor3-08}). In particular, this first result contributes an answer to a problem posed by F.~Grunewald, B.~Kunyavski\u{\i} and E.~Plotkin in \cite{GKP}. These  authors present a version of Theorem~\ref{G_S} for general classes of groups with good hereditary properties \cite[Theorem 5.12]{GKP}, by means of the following concepts:
\begin{de}\textup{(Grunewald, Kunyavski\u{\i}, Plotkin, \cite[Definition 5.10]{GKP})} For a class $\cal X$ of groups and a group $G$, an element $g\in G$ is called {\it locally $\cal X$-radical} if
$\langle g^{\langle x\rangle}\rangle$ belongs to $\cal X$ for every $x\in G$; and the element $g\in G$  is called {\it globally $\cal X$-radical} if
$\langle g^G\rangle$ belongs to $\cal X$.
\end{de}

 For a subset $S$ and a subgroup   $X$ of a group $G$, we set {$\langle S^{X}\rangle=\langle s^x\mid s\in S,\ x\in X\rangle$, the smallest $X$-invariant subgroup of $G$ containing $S$. For $g\in G$, we write $\langle g^X\rangle$ for $\langle \{g\}^X\rangle$. When $\cal X$ is a Fitting class,  the property $\langle g^{\langle x\rangle}\rangle\in \cal X$ is equivalent to $g\in \langle g,x\rangle_{\cal X}$, the $\cal X$-radical of $\langle g,x\rangle$, as the property $\langle g^G\rangle\in \cal X$ is equivalent to $g\in G_{\cal X}$,  and these properties are useful in the problem of characterizing elements forming $G_{\cal X}$.  As mentioned in  \cite[Section~5.1]{GKP}, a main problem is to determine classes $\cal X$ for which locally and globally $\cal X$-radical elements coincide. Corollary~\ref{Cor4-08} gives a positive answer for the class $\mathcal N^k$ of finite soluble groups of nilpotent length at most $k\ge 1$.

 When the Fitting  class $\cal X$  is in addition closed under extensions and contains all cyclic groups,  the condition $\langle g^{\langle x\rangle}\rangle\in \cal X$ is equivalent to $\langle g,x\rangle\in \cal X$, but this is not the case for important classes of groups, as  the class $\cal N$ of finite nilpotent groups, or more generally $\mathcal N^k$, $k\ge 1$. In this situation the condition $\langle g,x\rangle\in \cal X$ for all elements $x\in G$ may well not be equivalent to $g\in G_{\cal X}$, but still of interest as shown in Corollary~\ref{Cor3-08} in relation with the hypercenter.\medskip

 We shall adhere to the notation used in \cite{D-H} and we refer
also to that book for the basic results on classes of groups. In
particular,
$\pi (G)$ denotes the set of all primes dividing the order of the group
$G$. Also $\mathcal A$ and $\mathcal S_{\pi}$, $\pi$ a set of primes, denote the classes of abelian groups and soluble $\pi$-groups, respectively. For the class of all finite $\pi$-groups, the residual of any group $X$  is denoted $O^{\pi}(X)$, and $O_\pi(X)$ stands for the corresponding radical of $X$. If $\cal F$ is a class of groups, then $\cal NF$ is the class of groups which are extensions of a nilpotent normal subgroup by a group in $\cal F$.
\smallskip

We gather next our main results. The first one extends to the universe of finite groups results for soluble groups in \cite[Theorem 3]{GHP08b}, \cite[Theorem 1, Proposition 1]{GHP08} and \cite[Theorem 1]{GHP11}.

\begin{teor}\label{Theorem} Let $G=AB$ be a finite group, $A,B\le G$. Then: \begin{enumerate}\item $A, B$ are $\cal NA$-connected if and only if $[A,B]\le F(G)$.
\item $A, B$ are $\mathcal N^2$-connected if and only if $AF(G)/F(G)$ and $BF(G)/F(G)$ are
$\mathcal{N}$-connected.
\item Let $\pi,\rho$ be arbitrary sets of primes. The following are equivalent:
\begin{enumerate}[(a)]\item $A,B$ are $\mathcal S_\pi\mathcal S_\rho$-connected.
\item $\pi(G)\subseteq \pi\cup \rho$, $[A,B]$ is soluble, $[A,O^\rho(B)]\in \mathcal S_\pi$, $[B,O^\rho(A)]\in \mathcal S_\pi$.
\item $\pi(G)\subseteq \pi\cup \rho$, $[A,B]\in \mathcal S_\pi\mathcal S_\rho$, $[A,O^\rho(B)]\in \mathcal S_\pi$, $[B,O^\rho(A)]\in \mathcal S_\pi$.
\end{enumerate}

\item Let $\mathcal{F}$ be a formation of soluble groups containing all
abelian groups. Assume that one of the following conditions holds:
\begin{enumerate}[(i)]\item One of
the factors $A,B$ is normally embedded in $G$ (in the sense of \cite[I. Definition (7.1)]{D-H}).
\item $A$ and $B$ have coprime indices in $G$.
 \end{enumerate} Then $A$ and $B$ are
$\mathcal{NF}$-connected if and only if
$AF(G)/F(G)$ and $BF(G)/F(G)$ are
$\mathcal{F}$-connected.\end{enumerate}
\end{teor}

\begin{rem} In Theorem~\ref{Theorem}~(3),~(b) and (c), $[A,O^\rho(B)]\in \mathcal S_\pi$ is equivalent to $[A,O^\rho(B)]\le O_\pi(G)$, and also to
$[A,O^\rho(B)]\le O_\pi (G_\mathcal S)$ as $[A,B]$ is soluble.\smallskip

 This is because $[A,O^\rho(B)]$ is subnormal in $G$ since $[A,O^\rho(B)]\unlhd \langle O^\rho(B)^A\rangle = \langle O^\rho(B)^{BA}\rangle= \langle O^\rho(B)^G\rangle\unlhd G=AB$. \end{rem}

As consequences of Theorem~\ref{Theorem} we derive Corollaries~\ref{Cor1-08},~\ref{Cor2-08},~\ref{Cor3-08},~\ref{Cor4-08}, and point out again that corresponding results for finite soluble groups were firstly obtained in \cite[Corollaries 1, 2, 3, 4]{GHP08}.

\begin{cor}\label{Cor1-08}
Let  the group $G=AB$ be the $\mathcal{N}^{2}$-connected product of
the subgroups $A$ and $B$, and let ${\mathcal F}$ be a class of groups.
\begin{enumerate}\item  Assume that
\begin{description}
\item(i) ${\mathcal F}$ is a {\sc q}-closed Fitting
class, or \item(ii) ${\mathcal F}$ is either a saturated formation or a formation containing
${\mathcal N}$.
\end{description}
Then $A,B\in {\mathcal F}$ implies
$G\in {\mathcal N}{\mathcal F}$.
\item Assume that
\begin{description}
\item(i) ${\mathcal F}$ is a {\sc q}-closed Fitting
class,  or \item(ii') ${\mathcal F}$ is either a saturated formation or a formation of soluble groups containing
${\mathcal N}$.
\end{description}
Then $G\in {\mathcal N}{\mathcal
F}$ implies $A,B\in {\mathcal N}{\mathcal F}$.
\end{enumerate}
\end{cor}

As a particular case of Corollary~\ref{Cor1-08} we state explicitly:
\begin{cor}\label{Cor2-08}
If the group $G=AB$ is the $\mathcal{N}^{2}$-connected product of
the $\pi$-separable subgroups $A$ and $B$ of $\pi$-length at most $l$, $\pi$ a set of primes,
then $G$ is $\pi$-separable of $\pi$-length at most $l+1.$
\end{cor}

\begin{rem} Easy examples show that the bound for the $\pi$-length of $G$ in Corollary~\ref{Cor2-08} is sharp. For instance, for any $l\ge 1$, consider a set of primes $\pi\neq \emptyset$ with $\pi'\neq \emptyset$, where $\pi'$ stands for the complement of $\pi$ in the set of all prime numbers, let $B$ be a $\pi$-separable group of $\pi$-length $l$ such that $O_\pi(B)=1$, let $p\in \pi$ and $A$ be a faithful module for $B$ over the field of $p$ elements. Let $G=[A]B$ be the corresponding semidirect product of $A$ with $B$. Then $A$ and $B$ are $\mathcal{N}^{2}$-connected and the $\pi$-length of $G$ is $l+1$.
\end{rem}

\begin{cor} \label{Cor3-08}
Let $G$ be a group, $g\in G$ and $k\geq 1$. Then $\langle
g,h\rangle \in \mathcal N^k$ for all $h\in G$ if and only if $g\in
Z_\infty (G\mbox{ mod }F_{k-1}(G))$.
\end{cor}

\begin{rem}\label{hypercenter}\begin{enumerate}\item For $k=1$, Corollary~\ref{Cor3-08} gives a  characterization of the hypercenter of a group. This particular case appears already in \cite[Corollary 3]{GHP08} as a direct consequence of Lemma~\ref{lem:nco}~(2) below, and was also observed by R.~Maier, as mentioned in \cite[Remark 5.5]{GKP}, and referred to \cite{RM}.

\item Assume that $g\in G$ such that $\langle g,x\rangle$ is soluble for all $x\in G$. Let $l$ be the highest nilpotent length of all these subgroups, so that $\langle g,x\rangle\in \mathcal N^l$ for all $x\in G$. By Corollary~\ref{Cor3-08} it follows that $g\in F_l(G)\le G_{\mathcal S}$. So that Corollary~\ref{Cor3-08} may be seen also as generalization of the characterization of the soluble radical in Theorem~\ref{G_S}.
\end{enumerate}\end{rem}

\begin{cor} \label{Cor4-08}
For a finite group $G$, an element $g\in G$ and $k\geq 1$, the
following statements are equivalent:
\begin{enumerate}
\item $g\in F_k(\langle g,h\rangle)$ for all $h\in G$, i.e. $g$ is locally $\mathcal N^k$-radical. \item
$\langle g,h\rangle \in \mathcal N^k$ for all $h\in F_k(G)$ and
$\langle g,h\rangle \in \mathcal N^{k+1}$ for all $h\in G$. \item
$g\in F_k(G)$, i.e. $g$ is globally $\mathcal N^k$-radical.
\end{enumerate}
\end{cor}

\begin{rem} 1. P. Flavell  proved the equivalence of Condition 1 and
Condition 3 for $k=2$ (\cite[Theorem A]{FL02}) and  for arbitrary $k$ and soluble groups (\cite[Theorem 2.1]{FL02}),  as mentioned also in \cite[Remark 3]{GHP08}.
\smallskip

2. For $k=1$, the Baer-Suzuki theorem states that $F(G)=\{ g\in G\mid \langle g,g^x\rangle \in \mathcal N \text{ for every }x\in G\}$. But for $k=2$, one can not conclude that $g\in F_2(G)$ whenever $\langle g,g^x\rangle \in \mathcal N^2$ for all $x\in G$, as pointed out by  Flavell \cite{FL02}.
\end{rem}

\begin{rem} As application of Theorem~\ref{Theorem}, the hypothesis of solubility can be also omitted in Corollary~4 and Propositions~3, 4 of \cite{GHP08b}, especially in relation with saturated formations $\mathcal F\subseteq \cal NA$, such as the class of supersoluble groups. Also an extension for finite groups of Corollary~1 of \cite{GHP11}, in relation with the above-mentioned nilpotent-like Fitting formations, can be stated.
\end{rem}

\section{Proofs of the main results}

\begin{lem} \label{lem:nco}\textup{(\cite[Proposition~1~(2),(8), Lemma~1, Proposition~3]{Ncon2})} Let the group
$G=AB$ be an $\mathcal N$-connected product of the subgroups $A$
and $B$. Then:
\begin{enumerate}
\item $A$ and $B$ are subnormal in $G$. \item $A\cap B\le Z_\infty (G)$, the hypercenter of $G$. \item If $\cal F$ is either  a saturated formation or a formation containing $\cal N$, and $A,B\in \cal F$, then $G\in \cal F$.
\item If $\cal F$ is either a saturated formation or a formation of soluble groups containing $\cal N$, and $G\in \cal F$, then $A,B\in \cal F$.
\end{enumerate}

\end{lem}

\begin{pro}\label{Prop1-08}\textup{\cite[Proposition 1]{GHP08}}
Let $\mathcal{F}$ be a formation of soluble groups containing all
abelian groups. Let $G$ be a soluble group such that $G=AB$ is the
$\mathcal{NF}$-connected product of the subgroups $A$ and $B$.
Assume that one of the following conditions holds:
\begin{description}\item(i) One of the factors $A,B$ is normally embedded in $G$.
\item(ii) $A$ and $B$ have coprime indices in $G$.
\item(iii) $A$ and $B$ are nilpotent.
\end{description}
 Then
\[G/F(G)=(AF(G)/F(G))(BF(G)/F(G))\] is
an $\mathcal{F}$-connected product of the two factors.
\end{pro}

\noindent
{\bf Proof.} Part (i) is \cite[Proposition 1]{GHP08}. Parts (ii) and (iii) follow with the same arguments, taking into account that in both cases the following fact holds: If in addition the group $G$ has a unique minimal normal subgroup $N=C_G(N)$,  then either $N\le A$ or $N\le B$. (In particular, for part (iii) apply \cite[Theorem 1]{G73}.)\qed

\begin{lem}\label{lem1} Let $G=AB$ be a finite group, $A,B\le G$, and suppose that $[A,B]$ is soluble. Let $a\in A$. If $B\langle a\rangle G_{\cal S}=G$, then $A\le G_{\cal S}$.
\end{lem}

\noindent
{\bf Proof.} Since $[A,B]\le G_{\cal S}$ and $B\langle a\rangle G_{\cal S}=G$, it follows that $\langle a\rangle G_{\cal S}\trianglelefteq G$. Therefore
$\langle a\rangle \le G_{\cal S}$ and $BG_{\cal S}=G$. For any $x\in A$, we have now that $B\langle x\rangle G_{\cal S}=G$, and so again $x\in G_{\cal S}$, which implies $A\le G_{\cal S}$.\qed

\begin{de}\label{DefF} We define a \emph{subset functor} $T$ to assign to each finite group $G$ a subset $T(G)$ of $G$ satisfying the following conditions:
\begin{enumerate}\item $T(H)=\alpha(T(G))$ for all group isomorphisms $\alpha:G\longrightarrow H$.
\item $T(G)\cap U\subseteq T(U)$ for all groups $G$ and all $U\le G$.
\end{enumerate}
\end{de}

\begin{notation}\label{notation} For any group $G$, let $P(A,B,G)$ be a property on $G,A,B$, where $A,B$ are subgroups of $G=AB$, which satisfies the following conditions:
\begin{enumerate}\item Whenever $\alpha:G\longrightarrow H$ is a group isomorphism, if $P(A,B,G)$ is true, then $P(\alpha(A),\alpha(B),H)$ is true.
\item Whenever $L,A,B$ are subgroups of a group $G=AB$, with $L$ of the form $L=(L\cap A)(L\cap B)=G_{\cal S}X\langle y\rangle$, $\{X,Y\}=\{A,B\}$, $y\in Y$,
if $P(A,B,G)$ is true, then $P(L\cap A, L\cap B,L)$ is true.
\end{enumerate}
\end{notation}

\begin{pro}\label{proposition} Let $S_1,S_2$ be subset functors according to Definition~\ref{DefF}, $\mathcal Y\supseteq \cal A$ be a formation and $\cal X$ be a Fitting class.

Suppose the following statement $(\star)$ holds for all finite soluble groups.
\smallskip

$(\star)\begin{cases}\text{If}\ P(A,B,G)\ \text{holds in the group}\ G=AB,\\ \text{then}\ \langle a,b\rangle ^{\cal Y}\le G_{\cal X}\ \text{for all}\ a\in S_1(A), b\in S_2(B).\end{cases}$
\smallskip

Then $(\star)$ holds for all finite groups $G=AB$ such that $[A,B]$ is soluble.
\end{pro}
\smallskip
\noindent

{\bf Proof.} Let $G=AB$ be a finite group, $A,B\le G$,  such that $[A,B]\le G_{\cal S}$, and assume that $P(A,B,G)$ holds.  Let $a\in S_1(A), b\in S_2(B)$. We aim to  prove that $\langle a,b\rangle ^{\cal Y}\le G_{\cal X}$. We argue by induction on $|G|$. Assume first that $G_{\cal S}\langle a\rangle B=G$. Then $A\le G_{\cal S}$ by Lemma~\ref{lem1}. If $G=G_{\cal S}\langle b\rangle$, then $G$ is soluble and the result follows. So we may assume that $G_{\cal S}\langle b\rangle=A(B\cap G_{\cal S}\langle b\rangle)<G$. Since $b\in S_2(B)\cap (B\cap G_{\cal S}\langle b\rangle)\subseteq S_2(B\cap G_{\cal S}\langle b\rangle)$, and $G_{\cal S}\langle b\rangle$ has the desired form in~\ref{notation}(2), by induction we have that  $\langle a,b\rangle ^{\cal Y}\le (G_{\cal S}\langle b\rangle)_{\cal X}\cap G_{\cal S}\le (G_{\cal S})_{\cal X}\le G_{\cal X}$, because $\langle a,b\rangle ^{\cal Y}\le \langle a,b\rangle '=[\langle a\rangle,\langle b\rangle]\le G_{\cal S}.$ If $G_{\cal S}\langle a\rangle B=(G_{\cal S}\langle a\rangle B\cap A)B<G$, the same argument, with $G_{\cal S}\langle a\rangle B$ playing the role of $G_{\cal S}\langle b\rangle$, proves that $\langle a,b\rangle ^{\cal Y}\le G_{\cal X}$. (Note that $G_{\cal S}\langle a\rangle B$ is a subgroup of $G$ because $[A,B]\le G_{\cal S}$.) The proposition is proved.\qed
\medskip

\begin{rem} As we will see, Proposition~\ref{proposition} provides the main tool to derive Theorem~\ref{Theorem} from Theorem~\ref{MT} and the corresponding previous results in the soluble universe. In Notation~\ref{notation}~(2), the additional restriction of subgroups  $L=(L\cap A)(L\cap B)$ to subgroups of the form $L=(L\cap A)(L\cap B)=G_{\cal S}X\langle y\rangle$,  $\{X,Y\}=\{A,B\}$, $y\in Y$, will be required only for the application to the proof of Part (4) of Theorem~\ref{Theorem}, as it is also the case of the following Lemma~\ref{lem:ne}. The present formulations of Notation~\ref{notation} and Proposition~\ref{proposition} unify the treatment of the different parts stated in Theorem~\ref{Theorem}.\end{rem}
 \smallskip

\begin{lem}\label{lem:ne}\begin{enumerate}\item  Assume that $A$ is a normally embedded subgroup of a group $G=NA$ where $N\unlhd G$. Let $a\in A$. Then $N\langle a\rangle \cap A$ is normally embedded in $N\langle a\rangle$.

\item Assume that $A$ is a normally embedded subgroup of a group $G=AB$, $B\le G$, such that  $[A,B]$ is soluble. Then $G_{\cal S}B\langle a\rangle \cap A$ is normally embedded in $G_{\cal S}B\langle a\rangle$, for any $a\in A$.
    \item Let $G=AB$ a group such that $A$ and $B$ are subgroups of coprime indices in $G$, and $[A,B]$ is soluble. Then $G_{\cal S}B\langle a\rangle \cap A$ and $B$ have  coprime indices in  $G_{\cal S}B\langle a\rangle$, for any $a\in A$. \end{enumerate}
\end{lem}

\noindent
{\bf Proof.} 1. Let $p\in \pi (A)$. We consider $a=a_pa_{p'}$, where $a_p,\ a_{p'}$ denote the $p$-part and the $p'$-part of $a$, respectively. Let $M_p\in \Syl_p(A)$ such that $a_p\in M_p$. By the hypothesis, there exists $M\unlhd G$ such that $M_p\in \Syl_p(M)$.

We claim that $(N\cap M_p)\langle a_p\rangle \in \Syl_p (N\langle a\rangle \cap A)$. Since $N\cap A\trianglelefteq A$, we have that $N\cap M_p\in \Syl_p(N\cap A)$. Consequently, $(N\cap M_p)\langle a_p\rangle \in \Syl_p((N\cap A)\langle a_p\rangle)$. Since $N\langle a\rangle \cap A=(N\cap A)\langle a\rangle=(N\cap A)\langle a_p\rangle\langle a_{p'}\rangle$, the claim follows easily.

We prove next that $(N\cap M_p)\langle a_p\rangle \in \Syl_p (N\langle a\rangle \cap M)$. Since $N\langle a\rangle \cap M\trianglelefteq N\langle a\rangle$, it will follow that $N\langle a\rangle \cap A$ is normally embedded in $N\langle a\rangle$, which will conclude the proof.

We notice that $N\langle a\rangle \cap M=(N\langle a_{p'}\rangle \cap M)\langle a_p\rangle$, so that it is enough to prove that
$N\cap M_p\in \Syl_p(N\langle a_{p'}\rangle \cap M)$.

Again $N\cap M\trianglelefteq M$ implies that $N\cap M_p\in \Syl_p(N\cap M)$. Let $N_p\in Syl_p(N)$ such that $N\cap M_p=N_p\cap M_p$.  Then
$N\cap M_p=N_p\cap M_p\le N_p\cap M\in \Syl_p(N\cap M)$ because also $N\cap M\trianglelefteq N$. Consequently, $N\cap M_p=N_p\cap M_p= N_p\cap M$.

On the other hand, $N_p\in \Syl_p(N\langle a_{p'}\rangle)$ and $N\langle a_{p'}\rangle\cap M\trianglelefteq N\langle a_{p'}\rangle$, which implies that
$N_p\cap M\in \Syl_p(N\langle a_{p'}\rangle \cap M)$, and we are done.
\smallskip

2. Since $[A,B]\le G_{\cal S}$, we have that $BG_{\cal S}\trianglelefteq G=AB=BG_{\cal S}A$. The result follows now from part~1.
\smallskip

3. Set $N=BG_{\cal S}\trianglelefteq G=AB$, as before. Notice that
$$|N\langle a\rangle:N\langle a\rangle\cap A|=|N\langle a\rangle A:A|=|NA:A|=|G:A|.$$
Then $\gcd\,(|N\langle a\rangle :B|, |N\langle a\rangle :N\langle a\rangle\cap A|)\ |\ \gcd\,(|G:B|,|G:A|)=1$, and the result follows.\qed

\medskip

\noindent
{\bf Proof of Theorem~\ref{Theorem}.}
\smallskip

1. Apply Proposition~\ref{proposition} with $P(A,B,G)$ being $A,B$ $\mathcal{NA}$-connected, $S_1(G)=S_2(G)=G$ for all groups $G$, $\mathcal Y=\mathcal A$, $\mathcal X=\mathcal N$, Theorem~\ref{MT} and \cite[Theorem 3]{GHP08b}.
\smallskip

2. Apply Proposition~\ref{proposition} with $P(A,B,G)$ being $A,B$ $\mathcal{N}^2$-connected, $S_1(G)=S_2(G)=G$ for all groups $G$, $\mathcal Y=\mathcal N$, $\mathcal X=\mathcal N$, Theorem~\ref{MT} and \cite[Theorem 1]{GHP08}.
\smallskip

3. $(a) \Longrightarrow (b)$:\\
Apply Proposition~\ref{proposition} with $P(A,B,G)$ being $A,B$ $\mathcal{S_\pi S_\rho}$-connected, $S_1(G)=G$, $S_2(G)=\{g\in G\mid g\text{ is a }\rho'\text{-element}\}$,  for all groups $G$, $\mathcal Y=\mathcal A$, $\mathcal X=\mathcal \mathcal S_{\pi}$, Theorem~\ref{MT} and \cite[Theorem 1]{GHP11}.
\smallskip

\noindent
$(b) \Longrightarrow (c)$:\\
With the notation of Proposition~\ref{proposition}, let $P(A,B,G)$ be defined as follows:
\begin{itemize}\item $\pi(G)\subseteq \pi\cup \rho,\ [A,O^{\rho}(B)]\in \mathcal S_\pi, [B,O^{\rho}(A)]\in \mathcal S_\pi.$
\end{itemize}
In addition set $S_1(G)=S_2(G)=G$,   for all groups $G$, $\mathcal Y=\mathcal A$, $\mathcal X=\mathcal \mathcal S_{\pi}\mathcal S_{\rho}$.
\smallskip

We notice that in this case the condition $\langle a,b\rangle^{\cal Y}=\langle a,b\rangle '=[\langle a\rangle, \langle b\rangle]\le G_{\cal X}=G_{\mathcal S_\pi \mathcal S_\rho}$ for all $a\in A$ and $b\in B$, is equivalent to $[A,B]\in \mathcal S_\pi \mathcal S_\rho$.
\smallskip

We prove next that whenever $G=AB\in \cal S$, $\pi(G)\subseteq \pi\cup \rho$, $A,B\le G$, $[A,O^\rho(B)]\in \mathcal S_\pi$, $[B,O^\rho(A)]\in \mathcal S_\pi$, then $[A,B]\in \mathcal S_\pi \mathcal S_\rho$.
\smallskip

For such a group $G=AB$, we argue as in the proof of \cite[Theorem 1, (b)$\Rightarrow$(a)]{GHP11} and consider $A=O^\rho (A)A_\rho$  and $B=O^\rho (B)B_\rho$, where $A_\rho$ and $B_\rho$ are Hall $\rho$-subgroups of $A$ and $B$, respectively, such that $A_\rho B_\rho$ is a Hall $\rho$-subgroup of $G$. Then: $$[A,B]O_\pi(G)=[A, O^\rho (B)B_\rho]O_\pi(G)=[A, B_\rho]O_\pi(G)=[A_\rho, B_\rho]O_\pi(G).$$ Since $A_\rho B_\rho\in \mathcal S_\rho$, it follows that $[A,B]\in \mathcal S_\pi\mathcal S_\rho$.
\smallskip

We can apply now Proposition~\ref{proposition} to deduce that (b) implies (c).
\smallskip

\noindent
$(c) \Longrightarrow (a)$:\\
Apply Proposition~\ref{proposition} with $P(A,B,G)$ being:
\begin{itemize}\item $\pi(G)\subseteq \pi\cup \rho,\ [A,O^{\rho}(B)]\in \mathcal S_\pi, [B,O^{\rho}(A)]\in \mathcal S_\pi,$
\end{itemize}
$S_1(G)=S_2(G)=G$,   for all groups $G$, $\mathcal Y=\mathcal S_{\pi}\mathcal S_{\rho}\mathcal S_{(\pi\cup \rho)'}$, $\mathcal X=(1)$, and \cite[Theorem 1]{GHP11}.
\smallskip

4. Apply Proposition~\ref{proposition} with $P(A,B,G)$ being $A,B$ $\mathcal{NF}$-connected, either $A$ or $B$ normally embedded in $G$ for the case (i), or $A$ and $B$ of coprime indices in $G$ for the case (ii), $S_1(G)=S_2(G)=G$ for all groups $G$, $\mathcal Y=\mathcal F$, $\mathcal X=\mathcal N$, Lemma~\ref{lem:ne}, Theorem~\ref{MT} and Proposition~\ref{Prop1-08}.\qed
\medskip

\noindent
{\bf Proof of Corollary~\ref{Cor1-08}.} If $G=AB$ is an
${\mathcal N^2}$-connected product of subgroups $A$ and
$B$, then $G/F(G)$ is the $\mathcal{N}$-connected product of the
subgroups $AF(G)/F(G)$ and $BF(G)/F(G)$ by Theorem~\ref{Theorem}~(2). The result follows now from Lemma~\ref{lem:nco}.\qed
\medskip

\noindent
{\bf Proof of Corollary~\ref{Cor3-08}.} Mimic the proof of \cite[Corollary 3]{GHP08} by applying Theorem~\ref{Theorem}~(4)(i).\qed
\medskip

\noindent
{\bf Proof of Corollary~\ref{Cor4-08}.} Mimic the proof of \cite[Corollary 4]{GHP08} by using now Corollary~\ref{Cor3-08}.\qed

\bigskip

\noindent
{\bf Acknowledgments.} Research supported by Proyectos PROMETEO/2017/ 057 from the Generalitat Valenciana (Valencian Community, Spain), and PGC2018-096872-B-I00 from the Ministerio de Ciencia, Innovaci\'on y Universidades, Spain, and FEDER, European
Union; and third author also by Project VIP-008 of Yaroslavl P. Demidov State University.

\bigskip

\noindent
\footnotesize{M.~P. G\'{A}LLEGO}\\
\footnotesize{Departamento de Matem\'aticas, Universidad de Zaragoza,}\\
\footnotesize{Edificio Matem\'aticas, Ciudad Universitaria,
50009 Zaragoza, Spain}\\
 \footnotesize{E-mail: pgallego@unizar.es}\\
\\
\footnotesize{P. HAUCK}\\
\footnotesize{Fachbereich Informatik,
Universit\"{a}t T\"{u}bingen,}\\
\footnotesize{Sand 13, 72076 T\"{u}bingen, Germany}\\
\footnotesize{E-mail: peter.hauck@uni-tuebingen.de}\\
 \\
 \footnotesize{L. S. KAZARIN}\\
\footnotesize{Department of Mathematics, Yaroslavl P. Demidov State University}\\
 \footnotesize{Sovetskaya Str 14, 150014 Yaroslavl, Russia}\\
 \footnotesize{E-mail: Kazarin@uniyar.ac.ru}\\
 \\
 \footnotesize{A. MART\'{I}NEZ-PASTOR}\\
 \footnotesize{Instituto Universitario de Matem\'{a}tica Pura y  Aplicada IUMPA}\\
\footnotesize{Universitat Polit\`{e}cnica de Val\`{e}ncia,
 Camino de Vera, s/n,  46022 Valencia, Spain}\\
 \footnotesize{E-mail: anamarti@mat.upv.es}\\
 \\
 \footnotesize{and M.~D. P\'{E}REZ-RAMOS}\\
\footnotesize{Departament de Matem\`{a}tiques, Universitat de Val\`{e}ncia,}\\
 \footnotesize{C/ Doctor Moliner 50, 46100 Burjassot
({Val\`{e}ncia}), Spain}\\
\footnotesize{E-mail: Dolores.Perez@uv.es}

\end{document}